\documentclass[sn-mathphys,Numbered]{sn-jnl}

\usepackage{graphicx}%
\usepackage{multirow}%
\usepackage{amsmath,amssymb,amsfonts}%
\usepackage{amsthm}%
\usepackage{mathrsfs}%
\usepackage[title]{appendix}%
\usepackage[table,dvipsnames]{xcolor}%
\usepackage{textcomp}%
\usepackage{manyfoot}%
\usepackage{booktabs}%
\usepackage{algorithm}%
\usepackage{algorithmicx}%
\usepackage{algpseudocode}%
\usepackage{listings}%

\usepackage{cleveref}
\usepackage{caption}
\usepackage{subcaption}
\usepackage{siunitx}
\usepackage{tabularx}
\usepackage{comment}
\usepackage{arydshln}

\usepackage{tikz}
\usetikzlibrary{fit}
\usetikzlibrary{calc}
\pgfdeclarelayer{background}
\pgfsetlayers{background,main}
\colorlet{inbox}{lightgray!20}
\colorlet{outbox}{lightgray!50}

\newcommand{\vp}{\varphi}

\newcommand{\R}{\mathbb{R}}
\newcommand{\Rp}{\mathbb{R}_{\ge 0}}
\newcommand{\N}{\mathbb{N}}

\newcommand{\uFC}{u_{\rm fb}}

\newcommand{\sgn}{\rm sign}

\newcommand{\red}[1]{\textcolor{red}{#1}}
\newcommand{\ljl}[1]{\textcolor{teal}{#1}}

\graphicspath{{Figures}}

\theoremstyle{thmstyleone}%
\theoremstyle{thmstyletwo}%

\theoremstyle{thmstylethree}%

\begin{document}

\title[Experimental validation for the combination of funnel control with a feedforward control strategy]{Experimental validation for the combination of funnel control with a feedforward control strategy}


\author[1]{\fnm{Svenja} \sur{Dr\"ucker}}\email{svenja.druecker@tuhh.de}

\author[2]{\fnm{Lukas} \sur{Lanza}}\email{lukas.lanza@tu-ilmenau.de}

\author[3]{\fnm{Thomas} \sur{Berger}}\email{thomas.berger@math.upb.de}

\author[2]{\fnm{Timo} \sur{Reis}}\email{timo.reis@tu-ilmenau.de}

\author*[1]{\fnm{Robert} \sur{Seifried}}\email{robert.seifried@tuhh.de}

\affil*[1]{\orgdiv{Institute of Mechanics and Ocean Engineering}, \orgname{Hamburg University of Technology}, \orgaddress{\street{Ei\ss endorfer Str. 42}, \postcode{21073}, \city{Hamburg}, \country{Germany}}}

\affil[2]{\orgdiv{Institute of Mathematics}, \orgname{Technische Universität Ilmenau}, \orgaddress{\street{Ehrenbergstraße 29}, \postcode{98693}, \city{Ilmenau}, \country{Germany}}}

\affil[3]{\orgdiv{Institute of Mathematics}, \orgname{Paderborn University}, \orgaddress{\street{Warburger Str. 100}, \postcode{33098}, \city{Paderborn},  \country{Germany}}}


\abstract{Current engineering design trends, such as light-weight machines and human-machine-interaction, often lead to underactuated systems. Output trajectory tracking of such systems is a challenging control problem. Here, we use a two-design-degree of freedom control approach by combining funnel feedback control with feedforward control based on servo-constraints. We present experimental results to verify the approach and demonstrate that the addition of a feedforward controller mitigates drawbacks of the funnel controller. We also present new experimental results for the real-time implementation of a feedforward controller based on servo-constraints on a minimum phase system.   }


\maketitle

\section{Introduction}
Current engineering trends include the development of light-weight machines, many types of cable-driven manipulators and flexible joint robots. These systems have more degrees of freedom than independent control inputs and are called underactuated systems~\cite{LiuYu13,Seifried14}. Underactuation occurs naturally in the design of many mechanical systems. With the development of new mechanical designs and functionalities, control strategies must be designed accordingly in order to meet accuracy requirements. For underactuated systems, it is usually not possible to control all degrees of freedom independently and therefore, output trajectory tracking control of such systems is a challenging problem. 
In this study, we follow the popular two degree of freedom control methodology~\cite{Skogestad04} in order to solve the trajectory tracking control problem, since it can be applied efficiently to underactuated systems. In particular, the considered control strategy combines funnel control~\cite{ilchmann2002tracking,berger2021tracking} (as the feedback component) with feedforward control based on servo-constraints~\cite{BlajKolo04,Druecker22}. 
This combination has been shown in simulation to be a successful control strategy for underactuated multibody systems, cf.~\cite{berger2019combined,berger2021tracking}. \\

\noindent
In this article we present results in order to experimentally validate the two-design-degree of freedom control strategy applied to a torsional oscillator.
The oscillator consists of two flywheels connected by a rod, introduced in \Cref{sec:exp} and depicted in \Cref{fig:exp}.
This system can be considered as a simplified model of a drive-train, cf.~\cite{Pham19}. Therefore, experimentally validating control strategies with this setup is an important step towards application of the proposed controller to real systems such as drive-trains.

\subsection{Background and motivation}
In the feedforward branch, we aim to cancel out all known dynamics by using an inverse model as feedforward controller. Classical approaches, such as the Byrnes-Isidori normal form~\cite{isidori95}, are often burdensome to derive for complex underactuated multibody systems due to the involved algebraic manipulations of the system dynamics. In contrast, the servo-constraints approach is an efficient strategy to compute the inverse model of underactuated multibody systems~\cite{BlajKolo04,Druecker22}. For this purpose, the equations of motion are appended by so-called servo-constraints, which constrain the system output to a prescribed trajectory. The resulting set of differential-algebraic equations (DAEs) describes the inverse model. The solution of the DAE problem directly includes the control inputs which move the (nominal) system on the prescribed trajectory. The inverse model DAEs can be solved in real-time~\cite{Otto2018}. In the context of servo-constraints, the DAEs are often solved using the implicit Euler scheme~\cite{BlajKolo04,SeifriedBlajer13}, but there also exist results for higher order integration schemes, such as backwards differentiation formulas~\cite{FumagalliMasaratiEtAl10,MobergHanssen10}. This methodology has shown to be an efficient control strategy for complex underactuated multibody systems. Typical application examples lie in the class of differentially flat systems, such as cranes~\cite{AltmannBetschEtAl16,BetschEtAl09,Otto2018} or mass-spring chains~\cite{AltmannHeiland17,FumagalliMasaratiEtAl10}. However, the method is also applicable to minimum phase systems~\cite{SeifriedBlajer13} or non-minimum phase systems. Non-minimum phase systems involve a stable inversion, which is reformulated for the servo-constraints framework in~\cite{BruelsEtAl13} and applied to flexible manipulators in~\cite{BastosEtAl17,DrueckerSeifried23,StroehleBetsch22}. 
\newline
While there exist many theoretical results, only few experimental studies are documented for the application of servo-constraints. An experimental study is presented in \cite{BlajerDziewieckiEtAl09} for a mass-spring system and in \cite{Bestle05} for a small scale crane system. An experimental study of the real-time capabilities of the scheme is presented in \cite{Otto2018} for the differentially flat crane system. There are even less experimental results for the real-time application of servo-constraints on minimum phase systems. For example in \cite{MorlockEtAl21}, the method is applied to a flexible manipulator, for which stable internal dynamics is obtained by adding a counter-weight as well as using output relocation.

Of course, the nominal model which forms the basis for the inversion can never capture all physical effects. There always remain some uncertainties, which are for example classified into five categories in~\cite{ThrunBurgardEtAl06}. Therefore, there will be a tracking error when the feedforward strategy is applied to the real system. This tracking error has to be compensated by a suitable feedback strategy. 

As the feedback component in the present work we apply the so-called \emph{funnel controller}, first proposed in the seminal work~\cite{ilchmann2002tracking}. It is a high-gain feedback control strategy, which admits the following features.
First, it achieves that the output $y(t)$ of a system tracks a given reference signal~$y_{\rm ref}(t)$ with prescribed accuracy.
This means, that the tracking error $e(t) := y(t) - y_{\rm ref}(t)$ is guaranteed to evolve within given (possibly time-varying) bounds, i.e., $\| e(t) \| < \psi(t)$ for all $t \ge 0$, where $\psi$ is the error tolerance given by the engineer. 
Second, the controller is model-free in the sense that the feedback law requires only the instantaneous values of the error signal $e(t)$; in particular, no system parameters are utilized (for systems with higher relative degree the feedback law involves the higher order derivatives of the error signal, but no system parameters).
Funnel control proved to be a powerful tool for tracking problems in various applications. 
To name but a few, it has been investigated in control of industrial servo-systems~\cite{Hack17}, underactuated multibody systems~\cite{berger2021tracking,berger2019combined}, electrical circuits~\cite{berger2014zero,senfelds2014electrical}, peak inspiratory pressure~\cite{pomprapa2015periodic}, and adaptive cruise control~\cite{berger2020cruise}.
Moreover, even control of infinite-dimensional systems has been investigated. 
For instance, a boundary controlled heat equation in~\cite{reis2015funnel}, performing reference tracking of a moving water tank was studied in~\cite{berger2022funnel}, and defibrillation processes of the human heart was considered in~\cite{berger2021funnel}.
For a more comprehensive review and further applications see the recent work~\cite{BergIlch21}.
Moreover, in that article a funnel controller was designed, which achieves asymptotic exact tracking for unknown nonlinear systems with arbitrary relative degree.
The aspect of exact tracking has been further considered in~\cite{lanza2022exact}, where exact output tracking in predefined finite time is achieved with funnel control.
Since in real applications system data is only available at discrete sampling times, a feedback controller for sampled-data systems was developed in~\cite{lanza2023sampled}, which achieves tracking with predefined (possibly time-varying) accuracy.
Recently, the idea of funnel control was utilized to design a Model Predictive Control scheme~\cite{berger2020learning,denn2022funnel}, which involves a particular stage cost and thus achieves reference tracking with prescribed error tolerance with superior control performance compared to pure feedback controllers.
This MPC scheme was further equipped with a feedback loop in~\cite{berger2023robust} to obtain a robust controller, and a learning scheme was added in~\cite{lanza2023learning} to improve the underlying model during the control phase.

Experimental results involving funnel control are presented in~\cite{hackl2011funnel,hackl2012bang,hackl2013funnel,wang2016extended}.
In~\cite{hackl2011funnel} a funnel controller was combined with a PI controller to perform speed control of an electrical drive, where the measurement data is noisy.
In~\cite{hackl2012bang} the Bang-Bang funnel controller, proposed in~\cite{liberzon2010bang}, was extended with a neutral mode, and its functioning was experimentally verified with position control of an electrical drive.
In~\cite{hackl2013funnel} position and speed control of an electrical drive using a saturated funnel controller for systems with relative degree two was successfully tested in the laboratory.
In~\cite{wang2016extended} the application of a funnel controller combined with an extended state observer was experimentally verified for a permanent-magnet synchronous motor.

\subsection{Problem statement}

The objective of this contribution is to address the need for experimental validation of the two-component control strategy described above. In terms of the feedforward control strategy, experimental results for its real-time application are so far limited for a minimum phase system. In terms of the combined strategy, despite theoretical results~\cite{berger2019combined,berger2021tracking}, the lack of empirical evidence hinders the practical application of the proposed control strategy. Therefore, this research aims to provide experimental results that validate the efficiency and reliability of the control strategy under real-world conditions. By conducting comprehensive experiments and analyzing the obtained data, this study seeks to bridge the gap between theory and practice.

\subsection{Scope and contribution}

In our previous works~\cite{berger2019combined,berger2021tracking}, we analyze the described control strategy in a theoretical framework with supporting simulation results. In this contribution, we provide experimental data which support the theoretical results. We apply the control strategy described above to a torsional oscillator with two rotating flywheels. We present experimental results which validate the real-time application of servo-constraints for a minimum phase system. Further, we present experimental data to show that the proposed combined control strategy can compensate the known drawbacks of the funnel controller, in particular, peaky input signals cf.~\cite{denn2022funnel}.
In addition, we present insights about the discrete implementation of funnel control on a test bench. 
This aspect was also considered in the numerical example in~\cite{berger2020learning}.
Thus, we take another step towards the application of funnel control to real-world applications.  

\subsection{Organization of the paper}

The remainder of the paper is organized as follows. First, we briefly introduce the modeling approach based on multibody system dynamics and the considered control strategy in \Cref{sec:modeling}. Afterwards, the feedforward control strategy based on servo-constraints is introduced in \Cref{sec:ffw}, while the funnel controller is presented in \Cref{sec:fb}. In \Cref{sec:exp_setup}, we describe the experimental setup and analyze its dynamics and properties with respect to the controller design. The experimental results are shown in \Cref{sec:exp}. Finally, we conclude the paper with summarizing remarks and an outlook in \Cref{sec:conclusion}.

\section{Modeling and control strategy}\label{sec:modeling}

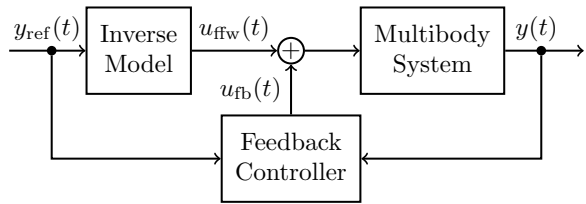
\begin{figure}[b]
  \centering
  \resizebox{0.6\textwidth}{!}{
\begin{tikzpicture}[thick,node distance = 12ex, box/.style={fill=white,rectangle, draw=black}, blackdot/.style={inner sep = 0, minimum size=3pt,shape=circle,fill,draw=black},plus/.style={fill=white,circle,inner sep = 0, minimum size=5pt,thick,draw},metabox/.style={inner sep = 3ex,rectangle,draw,dotted,fill=gray!20!white}]

    \node (sys)     [box,minimum size=8ex] {$\begin{array}{c} \text{Multibody}\\ \text{System}\end{array}$};
    \node (ufork)   [plus, left of = sys, xshift=-1ex] {$+$};
    \node (invsys)  [box,minimum size=8ex,left of = ufork,xshift=-2ex] {$\begin{array}{c} \text{Inverse}\\ \text{Model}\end{array}$};
    \node (yref)    [blackdot, left of = invsys, xshift=4ex] {};
    \node (yref1)   [minimum size=0pt, inner sep = 0pt, left of = yref, xshift=8ex] {};
    \node (fun)     [box, below of = ufork, yshift=2ex, xshift=0ex, minimum size=8ex]  {$\begin{array}{c} \text{Feedback}\\ \text{Controller}\end{array}$};
    \node (right1)  [blackdot,right of = sys,xshift=-2ex] {};
    \node (end)     [minimum size=0pt, inner sep = 0pt, right of = right1, xshift=-8ex] {};

    \draw[->] (ufork)   -- (sys)  {};  
    \draw[->]  (sys)     -- (end)  node[pos=0.4,above] {$y(t)$};
    \draw[->] (invsys)   -- (ufork)    node[midway,above] {$u_{\rm ffw}(t)$};
    \draw[->] (fun)   -- (ufork)    node[pos=0.5,left] {$u_{\rm fb}(t)$};
    \draw[->] (yref1) -- (invsys)    node[pos=0.5,above] {$y_{\rm ref}(t)$};
    \draw[->]  (yref)     |- (fun)  {};
    \draw[->]  (right1)     |- (fun)  {};

  \end{tikzpicture}
  }
\caption{Two design degree of freedom control approach for multibody systems.}
\label{fig:blockdiagram}
\end{figure}

We model the considered mechanical system using concepts from multibody dynamics~\cite{SchiEber14}. Here, we consider holonomic systems with $n$ degrees of freedom in minimal coordinates which are described by generalized coordinates~$q$. The equations of motion are given by
\begin{equation}\label{eq:MBS}
\begin{aligned}
    \dot q(t) &= v(t),\\
    M(q(t)) \dot v(t) &= f\big(q(t),v(t)\big) 
      + B(q(t))\, u(t),\\
    y(t) &=h\big(q(t),v(t)\big)
\end{aligned}
\end{equation}
with
\begin{itemize}
  \item the generalized coordinates~$q: I\to \R^n$ and generalized velocities~$v: I\to \R^n$ , where $I\subseteq \R_{\ge 0}$ is some interval,
  \item the generalized mass matrix $M:\R^n\to\R^{n\times n}$,
  \item the generalized forces $f:\R^n\times \R^n\to\R^n$,
  \item the input distribution matrix $B:\R^n\to\R^{n\times m}$,
  \item the output measurement function $h:\R^n\times\R^n\to\R^m$.
\end{itemize}
The functions $u:\Rp\to\R^m$ describe the inputs which act on the multibody system~\eqref{eq:MBS}. For underactuated systems, we have $m<n$ and it is not possible to control all degrees of freedom independently. The functions $y:\Rp\to\R^m$ are the outputs associated with the multibody system~\eqref{eq:MBS}. The functions in~\eqref{eq:MBS} are assumed to be sufficiently smooth and the mass matrix is assumed to be pointwise symmetric positive definite, i.e.,
\begin{equation}\label{eq:assumptions}
  \forall\,q\in\R^n:\quad M(q) = M(q)^\top>0 \,.
\end{equation}
The control objective is that the system output~$y$ tracks a prescribed trajectory $y_{\rm ref}:\R_{\ge 0}\to\R^m$. We follow the popular two-design-degree of freedom control structure, which combines a feedforward controller with a feedback controller, see e.g.~\cite{Skogestad04}. A block diagram of the methodology is shown in \Cref{fig:blockdiagram}. While the feedforward controller is responsible for steering the system output to the reference trajectory and to account for most of the motion, the feedback controller rejects disturbances and accounts for modeling errors. As in our preliminary simulation-based work~\cite{berger2019combined,berger2021tracking}, we choose a feedforward control strategy based on the method of servo-constraints and a funnel controller as the feedback component. Both concepts are briefly introduced in the following.

\section{Controller components} \label{Sec:ControllerComponents}
In this section we introduce the two controller components.
First, the feedforward component based on the concept of \emph{servo-constraints} is explained. Further, we discuss its implementation in the current context.
Second, the feedback component \emph{funnel control} is introduced.

\subsection{Feedforward control based on servo-constraints}\label{sec:ffw}
The feedforward controller is computed in terms of an inverse model of the system in order to cancel out all known dynamics. The method of servo-constraints provides an efficient and very general approach for complex underactuated multibody systems~\cite{BlajKolo04,Druecker22}.
Motivated from modeling classical mechanical constraints, such as joints, the equations of motion~\eqref{eq:MBS} are appended by $m$ servo-constraints
\begin{equation}\label{eq:servo-constraints}
    h(q(t),v(t)) - y_{\rm ref}(t) = 0,
\end{equation}
which enforce the output to stay on the prescribed trajectory $y_{\rm ref}$
. This results in the DAEs
\begin{equation}\label{eq:servo-constraintsDAE}
\begin{aligned}
    \dot q(t) &= v(t),\\
    M(q(t)) \dot v(t) &= f\big(q(t),v(t)\big) 
    + B(q(t))\, u(t),\\
    0 &=   h\big(q(t),v(t)\big) - y_{\rm ref}(t),
\end{aligned}
\end{equation}
which have to be solved numerically for the coordinates $q$ and $v$ as well as the input~$u$. For minimum phase systems, the DAEs~\eqref{eq:servo-constraintsDAE} can be integrated forward in time. For non-minimum phase systems, integration forward in time is not possible, since the states of the internal dynamics can become unbounded. In this case, a boundary value problem must be formulated in order to compute a bounded solution for the internal dynamics~\cite{BruelsEtAl13,ChenPaden96}. 
In the present article, we consider minimum phase systems and solve an initial value problem for~(\ref{eq:servo-constraintsDAE}).
The initial values $q(0), v(0), u(0)$ for \eqref{eq:servo-constraintsDAE} must be chosen so that they are consistent and the desired trajectory~$y_{\rm ref}$ must be compatible with the possible motion of the system, i.e., it is required that a solution of~\eqref{eq:servo-constraintsDAE} exists on $\R_{\ge 0}$. 

In order to compute the feedforward control input, we solve the inverse model DAEs~(\ref{eq:servo-constraintsDAE}) using the implicit Euler scheme
\begin{equation}\label{eq:implEul}
\begin{aligned}
    \hat q(t_{n+1}) &= \hat q(t_n) + \Delta t \, \hat v(t_{n+1}),\\
    \hat v(t_{n+1}) &= \hat v(t_{n}) + \Delta t \, M(\hat q(t_{n+1}))^{-1} \left( f\big(\hat q(t_{n+1}),\hat v(t_{n+1})\big) + B(\hat q(t_{n+1}))\, \hat u(t_{n+1}) \right) ,\\
    0 &=   h\big(\hat q(t_{n+1}),\hat v(t_{n+1})\big) - y_{\rm ref}(t_{n+1})
\end{aligned}
\end{equation}
with step size $\Delta t$ and the numerical approximations $\hat q(t_{n})$, $\hat v(t_{n})$ of the solution at time $t_{n}$. 
This integration scheme is a common choice in the context of servo-constraints~\cite{BlajKolo04,SeifriedBlajer13} because of its simplicity and its real-time capabilities.  Starting with initial values~$\hat q(t_{0}),\hat v(t_{0}),\hat u(t_{0})$, the nonlinear set of equations~(\ref{eq:implEul}) is solved for the solution~$\hat q(t_{n+1}),\hat v(t_{n+1}),\hat u(t_{n+1})$ at the next time instance in each control loop iteration. For this purpose Newton's method with a maximum of 10 iterations is applied.  

The numerical solution contains the trajectory for $u$, which is directly used as feedforward control, $u_{\rm ffw}:=u$. For a theoretical ideal model with exact parameters, this feedforward control ensures exact tracking with $y=y_{\rm ref}$. However, in practical applications model mismatches, disturbances and parameter uncertainties result in tracking errors. Reducing the latter effects is the purpose for adding a feedback controller, which is introduced in the next section. 

\subsection{Funnel control feedback component}\label{sec:fb}
\begin{figure}[h]
  \begin{center}
\begin{tikzpicture}[scale=0.35]
\tikzset{>=latex}
  \filldraw[color=gray!25] plot[smooth] coordinates {(0.15,4.7)(0.7,2.9)(4,0.4)(6,1.5)(9.5,0.4)(10,0.333)(10.01,0.331)(10.041,0.3) (10.041,-0.3)(10.01,-0.331)(10,-0.333)(9.5,-0.4)(6,-1.5)(4,-0.4)(0.7,-2.9)(0.15,-4.7)};
  \draw[thick] plot[smooth] coordinates {(0.15,4.7)(0.7,2.9)(4,0.4)(6,1.5)(9.5,0.4)(10,0.333)(10.01,0.331)(10.041,0.3)};
  \draw[thick] plot[smooth] coordinates {(10.041,-0.3)(10.01,-0.331)(10,-0.333)(9.5,-0.4)(6,-1.5)(4,-0.4)(0.7,-2.9)(0.15,-4.7)};
  \draw[thick,fill=lightgray] (0,0) ellipse (0.4 and 5);
  \draw[thick] (0,0) ellipse (0.1 and 0.333);
  \draw[thick,fill=gray!25] (10.041,0) ellipse (0.1 and 0.333);
  \draw[thick] plot[smooth] coordinates {(0,2)(2,1.1)(4,-0.1)(6,-0.7)(9,0.25)(10,0.15)};
  \draw[thick,->] (-2,0)--(12,0) node[right,above]{\normalsize$t$};
  \draw[thick,dashed](0,0.333)--(10,0.333);
  \draw[thick,dashed](0,-0.333)--(10,-0.333);
  \node [black] at (0,2) {\textbullet};
  \draw[->,thick](4,-3)node[right]{\normalsize$\inf\limits_{t \ge 0} \psi(t)$}--(2.5,-0.4);
  \draw[->,thick](3,3)node[right]{\normalsize$(0,e(0))$}--(0.07,2.07);
  \draw[->,thick](9,3)node[right]{\normalsize$\psi(t)$}--(7,1.4);
\end{tikzpicture}
\end{center}
 \caption{Evolution of the tracking error $e(t) = y(t)-y_{\rm ref}(t)$ within funnel boundary~$\psi(t)$; the figure is based on~\cite[Fig.~1]{BergLe18a}, edited for present purpose.}
 \label{Fig:FunnelPicture}
 \end{figure}
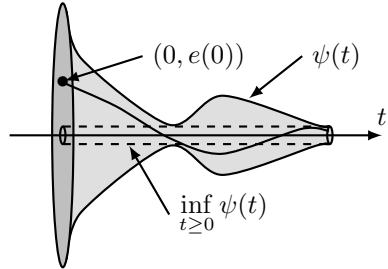

We briefly introduce the feedback-component of the control strategy, namely the \emph{funnel controller}, first proposed in~\cite{ilchmann2002tracking}. Although generalizations in several directions have been developed (see e.g.~\cite{BergIlch21} and the references therein), the version from~\cite{ilchmann2002tracking} for systems of relative degree one suffices for our purposes. For unknown nonlinear dynamical multi-input multi-output systems 
\begin{equation*}
    \begin{aligned}
        \dot x(t) &= f(x(t)) + g(x(t)) u(t), \quad x(0) = x^0 \in \R^n, \\
        y(t) &= h(x(t))
    \end{aligned}
\end{equation*}
with $h'(x) g(x)$ invertible for all $x \in \R^n$, i.e. relative degree 1, stable internal dynamics, and $g(x)$ sign definite ($v^\top g(x) v > 0$ or $v^\top g(x) v < 0$ for $v \in \R^m\setminus \{0\}$) for all~$x \in \R^n$,
the funnel controller achieves that the output~$y(t)$ follows a given reference~$y_{\rm ref}(t)$ with prescribed accuracy.
The latter means that the tracking error satisfies
\begin{equation}
  \forall \, t \ge 0 \, : \  \| y(t) - y_{\rm ref}(t) \| < \psi(t),
\end{equation}
where $\psi$ is a Lipschitz continuous and bounded function with $\inf_{s \ge 0} \psi(s) > 0$.
The situation is depicted in \Cref{Fig:FunnelPicture}.

The feedback law has the strikingly simple form
\begin{equation} \label{eq:uFC}
    \uFC (t) := -\frac{\psi(t)^2 (y(t) - y_{\rm ref}(t) )}{\psi(t)^2 - \|y(t) -y_{\rm ref}(t)\|^2} .
\end{equation}
Note that no system parameters are incorporated in the feedback law. 
Only the input/output dimension~$m \in \N$ is assumed to be known, and the structural assumptions that the system is minimum phase (stability of the internal dynamics), satisfies $h'(x) g(x) > 0$ for all $x \in \R^n$, and $g(\cdot)$ sign definite is assumed. 

The intuition behind the feedback law~\eqref{eq:uFC} is as follows.
Whenever the output tracking performs well, not much input action is required, i.e., a small tracking error $y(t) - y_{\rm ref}(t)$ results in a small signal~$\uFC(t)$.
If, however, the tracking error approaches the prescribed tolerance, i.e., if $\| y(t) - y_{\rm ref}(t) \| \to \psi(t)$, then the expression $1/(\psi(t)^2 - \| y(t) - y_{\rm ref}(t) \|^2 )$ grows rapidly and ``pushes the error away from the boundary~$\psi(t)$''.
Note that $1/(\psi(t)^2 - \| y(t) - y_{\rm ref}(t) \|^2 )$ exceeds any bounded value for $\| y(t) - y_{\rm ref}(t) \| \to \psi(t)$.
This fact is used in the feasibility proof~\cite{ilchmann2002tracking} to show that the closed-loop system has a global solution, and all signals remain bounded.

\section{Experimental setup and model}\label{sec:exp_setup}

\begin{figure}[h]
     \centering
     \begin{subfigure}[b]{0.65\textwidth}
         \centering
\begingroup%
  \makeatletter%
  \providecommand\color[2][]{%
    \errmessage{(Inkscape) Color is used for the text in Inkscape, but the package 'color.sty' is not loaded}%
    \renewcommand\color[2][]{}%
  }%
  \providecommand\transparent[1]{%
    \errmessage{(Inkscape) Transparency is used (non-zero) for the text in Inkscape, but the package 'transparent.sty' is not loaded}%
    \renewcommand\transparent[1]{}%
  }%
  \providecommand\rotatebox[2]{#2}%
  \newcommand*\fsize{\dimexpr\f@size pt\relax}%
  \newcommand*\lineheight[1]{\fontsize{\fsize}{#1\fsize}\selectfont}%
  \ifx\svgwidth\undefined%
    \setlength{\unitlength}{210.2510301bp}%
    \ifx\svgscale\undefined%
      \relax%
    \else%
      \setlength{\unitlength}{\unitlength * \real{\svgscale}}%
    \fi%
  \else%
    \setlength{\unitlength}{\svgwidth}%
  \fi%
  \global\let\svgwidth\undefined%
  \global\let\svgscale\undefined%
  \makeatother%
  \begin{picture}(1,0.52689682)%
    \lineheight{1}%
    \setlength\tabcolsep{0pt}%
    \put(0,0){\includegraphics[width=\unitlength,page=1]{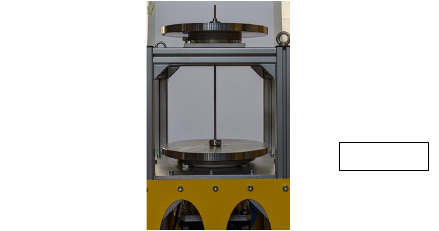}}%
    \put(0.78626486,0.15543231){\color[rgb]{0,0,0}\makebox(0,0)[lt]{\smash{\begin{tabular}[t]{l}actuator\end{tabular}}}}%
    \put(0,0){\includegraphics[width=\unitlength,page=2]{ExpSetup.pdf}}%
    \put(0.78632822,0.45302943){\color[rgb]{0,0,0}\makebox(0,0)[lt]{\smash{\begin{tabular}[t]{l}flywheels\end{tabular}}}}%
    \put(0,0){\includegraphics[width=\unitlength,page=3]{ExpSetup.pdf}}%
    \put(0.0168054,0.32503808){\color[rgb]{0,0,0}\makebox(0,0)[lt]{\smash{\begin{tabular}[t]{l}connecting \end{tabular}}}}%
    \put(0.02469803,0.27223678){\color[rgb]{0,0,0}\makebox(0,0)[lt]{\smash{\begin{tabular}[t]{l}shaft\end{tabular}}}}%
  \end{picture}%
\endgroup%

         \caption{Experimental setup.}
         \label{fig:exp}
    \end{subfigure}
     \hfill
     \begin{subfigure}[b]{0.34\textwidth}
         \centering
\begingroup%
  \makeatletter%
  \providecommand\color[2][]{%
    \errmessage{(Inkscape) Color is used for the text in Inkscape, but the package 'color.sty' is not loaded}%
    \renewcommand\color[2][]{}%
  }%
  \providecommand\transparent[1]{%
    \errmessage{(Inkscape) Transparency is used (non-zero) for the text in Inkscape, but the package 'transparent.sty' is not loaded}%
    \renewcommand\transparent[1]{}%
  }%
  \providecommand\rotatebox[2]{#2}%
  \newcommand*\fsize{\dimexpr\f@size pt\relax}%
  \newcommand*\lineheight[1]{\fontsize{\fsize}{#1\fsize}\selectfont}%
  \ifx\svgwidth\undefined%
    \setlength{\unitlength}{64.89781061bp}%
    \ifx\svgscale\undefined%
      \relax%
    \else%
      \setlength{\unitlength}{\unitlength * \real{\svgscale}}%
    \fi%
  \else%
    \setlength{\unitlength}{\svgwidth}%
  \fi%
  \global\let\svgwidth\undefined%
  \global\let\svgscale\undefined%
  \makeatother%
  \begin{picture}(1,1.4702817)%
    \lineheight{1}%
    \setlength\tabcolsep{0pt}%
    \put(0,0){\includegraphics[width=\unitlength,page=1]{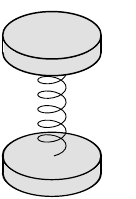}}%
    \put(0.61845336,0.50516223){\color[rgb]{0.25490196,0.46666667,0.24705882}\makebox(0,0)[lt]{\lineheight{0}\smash{\begin{tabular}[t]{l}$\varphi_1$\end{tabular}}}}%
    \put(0.60653282,1.4018778){\color[rgb]{0.25490196,0.46666667,0.24705882}\makebox(0,0)[lt]{\lineheight{0}\smash{\begin{tabular}[t]{l}$\varphi_2$\end{tabular}}}}%
    \put(0,0){\includegraphics[width=\unitlength,page=2]{MassChain_model.pdf}}%
    \put(0.81547261,0.27137697){\color[rgb]{0.63529412,0.10196078,0.10980392}\makebox(0,0)[lt]{\lineheight{0}\smash{\begin{tabular}[t]{l}$u$\end{tabular}}}}%
    \put(0,0){\includegraphics[width=\unitlength,page=3]{MassChain_model.pdf}}%
    \put(0.58245563,0.69245241){\color[rgb]{0,0,0}\makebox(0,0)[lt]{\lineheight{1.25}\smash{\begin{tabular}[t]{l}$k,d$\end{tabular}}}}%
  \end{picture}%
\endgroup%

         \caption{Multibody Model.}
         \label{fig:model}
     \end{subfigure}
     \hfill
        \caption{Experimental setup and multibody model.}
        \label{fig:ffw_vel}
\end{figure}

The control strategy described in \Cref{sec:modeling}, with components introduced in \Cref{Sec:ControllerComponents}, has so far been tested in simulations~\cite{berger2019combined,berger2021tracking}. In this contribution, we validate the control scheme experimentally on a torsional oscillator with two flywheels and a connecting rod. This is a simplified model of a drive train, e.g. such as~\cite{Pham19}. 
A picture of the experimental setup is shown in \Cref{fig:exp}. The lower flywheel is attached to the direct-drive actuator. The connecting shaft with $6\,$mm diameter connects the first and second flywheel. The system is modeled by two rigid bodies with inertia $I_1,\,I_2$, which are connected by a linear spring-damper combination with coefficients $k,\,d$, see \Cref{fig:model}. The rotation of the flywheels is described by the angle variables~$\varphi_1,\,\varphi_2$ and therefore, the generalized coordinates are chosen as $q=\begin{pmatrix} \varphi_1,  \varphi_2 \end{pmatrix}^\top$ and $v=\begin{pmatrix} \dot \varphi_1, \dot \varphi_2 \end{pmatrix}^\top$.  
The parameters of the model~(\ref{eq:sys}) are identified for the experimental setup shown in~\Cref{fig:exp} and are listed in \Cref{tab:param}. 

\begin{table}[h]
\begin{tabular}{cc}
\hline
& Nominal model \\
\hline
$I_1$ & \SI{0.136}{\kilogram\metre\squared}  \\
$I_2$ & \SI{0.12}{\kilogram\metre\squared}  \\
$k$& \SI{33.6}{\newton\metre} \\
$d$ & \SI{0.016}{\newton\metre\per\second} \\
\hline
\end{tabular}
\caption{Nominal model parameters.}
\label{tab:param}
\end{table}

According to the multibody dynamics approach described in~\Cref{sec:modeling}, the equations of motion are given by
\begin{equation} \label{eq:sys}
\begin{aligned}
    \begin{bmatrix} I_1 & 0 \\ 0 & I_2 \end{bmatrix}
    \begin{pmatrix} 
    \ddot \vp_1(t) \\ \ddot \vp_2(t)
    \end{pmatrix} &=
    \begin{bmatrix}
        -d & d \\ d & - d
    \end{bmatrix}
    \begin{pmatrix}
        \dot \vp_1(t) \\ \dot \vp_2(t)
            \end{pmatrix}
        +
        \begin{bmatrix}
            -k & k \\ k & -k
        \end{bmatrix}
        \begin{pmatrix}
            \vp_1(t) \\ \vp_2(t)
        \end{pmatrix}
        + \begin{bmatrix}
            F_{\rm f,1}(\dot \varphi_1(t))\\ 0
        \end{bmatrix}
        + \begin{bmatrix}
            1 \\ 0
        \end{bmatrix} u(t)
\end{aligned}
\end{equation}
with friction $F_{\rm f,1}$ acting on the first degree of freedom. 
The system input $u$ is given as a torque acting on the first flywheel. The constant model parameters~$I_1,\,I_2$ can be accurately measured and the parameters $\,k,\,d$ can be estimated with sufficient accuracy. However, the friction force $F_{\rm f,1}$ is more difficult to estimate. A simple approximation is given by assuming constant Coulomb friction, i.e., $F_{\rm f, 1}(s) = \tilde F \, \sgn(s)$, with $\tilde F > 0$, and $\sgn(\cdot)$ is the sign function.

We choose the angular velocity of the first flywheel as system output, which is to be tracked. Therefore, we have $y(t)=\dot \vp_1(t)$.
%
%
\newline

For analysis of the system, we introduce the relative motion $\Delta \vp := \vp_1 - \vp_2$. Thus, the rigid body motion is removed from the equations of motion~\eqref{eq:sys}, cf.~\cite{Pham19}.
We define the matrices
\begin{equation} \label{eq:reduced_system}
\begin{aligned}
   & M := \begin{bmatrix} 1&0&0\\ 0&I_1&0 \\ 0&0&I_2 \end{bmatrix}, \
    \tilde A := \begin{bmatrix} 0&1&-1\\-k&-d&d\\k&d&-d \end{bmatrix}, \
    \tilde B: = \begin{bmatrix}0\\1\\0\end{bmatrix}, \\
    & A : = M^{-1} \tilde A, \
    B := M^{-1} \tilde B,
    \end{aligned}
\end{equation}
and 
we set $F_1(s) := F_{\rm f,1}(s) / I_1$ and $F(s): = (0,F_1(s), 0)^\top$.
With $x:= (\Delta \vp, \dot \vp_1, \dot \vp_2)^\top$ the equations of motion for the reduced dynamics (the rigid body motion is removed) are given by
\begin{equation}
    \begin{aligned}
        \dot x(t) &= A x(t) + F(x(t)) + Bu(t) , \quad x(0) = x^0 \in \R^3, \\
        y(t) &= C x(t) = [0,1,0]x(t) = \dot \vp_1(t).
    \end{aligned}
\end{equation}
This dynamics describes the vibration around the rigid body motion. Since $C B = 1/I_1 \neq 0$, system~\eqref{eq:reduced_system} has 
relative degree~$r=1$.
In the next step, we decouple the internal dynamics from the equations of motion in order to analyze their stability.
Invoking the findings from~\cite{isidori95} and its extension for linear systems with bounded disturbances in~\cite[Ch.~3]{lanza2022contributions}, we may equivalently represent~\eqref{eq:reduced_system} as
\begin{equation} \label{eq:IO_reduced_system}
    \begin{aligned}
        \dot y(t) & = R y(t) + S \eta(t) + F_1(y(t)) + \Gamma u(t), \\
        \dot \eta(t) &= Q \eta(t) + P y(t) ,
    \end{aligned}
\end{equation}
where $\Gamma = CB = 1/I_1$, and
\begin{equation*}
    R = \frac{-d}{I_1}, \ S = \frac{1}{I_1} \begin{bmatrix} k & d \end{bmatrix}, \ 
    Q = \frac{1}{I_2} \begin{bmatrix} 0&I_2\\-k& -d \end{bmatrix}, \ 
    P = \frac{1}{I_2} \begin{bmatrix} -I_2 \\d \end{bmatrix}. 
\end{equation*}
In~\eqref{eq:IO_reduced_system}, $y = \dot \vp_1$ is the output and 
$\eta$ is the internal state given as $\eta = (-\Delta \vp, \dot \vp_2)^\top $ which cannot be seen directly in the output.
Since all eigenvalues of the matrix~$Q$ have negative real part, the internal dynamics of system~\eqref{eq:reduced_system} are bounded-input bounded-output stable, i.e., the system is minimum phase.
In particular, the internal dynamics are not directly effected by the roughly estimated Coulomb friction term.
Therefore, the control concepts described in~\Cref{Sec:ControllerComponents} are applicable to the model.

\section{Experimental results}\label{sec:exp}

\begin{figure}[h]
         \centering
         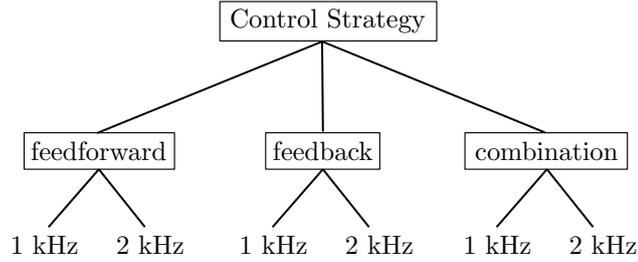
        \caption{Overview of the tested controllers.}
        \label{fig:overview}
\end{figure}

We performed a series of experiments on the experimental setup in order to validate the control scheme. An overview of the experiments is shown in \Cref{fig:overview}. The experimental setup offers two different control frequencies of $\SI{1}{kHz}$ and $\SI{2}{kHz}$. Advantages and disadvantages of the controllers are discussed for both frequencies. 
We first present measurements for the pure feedforward control and the pure feedback control, respectively. Afterwards, the combined control strategy is validated. We conclude this section by comparing all measurements using two different quantitative metrics.  

For all experiments, the desired trajectory is chosen as the polynomial
\begin{align*}
y_{\rm ref}(t) = \begin{cases}
y_{\rm 0}  , \quad &t < t_0 \\
y_{\rm 0} + \sigma(t)\left( y_{\rm f}-y_{\rm 0} \right) , \quad &t_0\leq t\leq t_{\rm f} \\
y_{\rm f}  , \quad &t> t_{\rm f}
\end{cases}
 \end{align*}
parameterized by the scalar parameter~$\sigma(t)$. The timing law of~$\sigma(t)$ is chosen as the polynomial
\begin{align*}
\sigma(t) = \,&\num{-3432}\left(\dfrac{t}{t_f}\right)^{15} +\num{25740}\left(\dfrac{t}{t_f}\right)^{14} - \num{83160}\left(\dfrac{t}{t_f}\right)^{13} +\num{150150}\left(\dfrac{t}{t_f}\right)^{12} \\ &- \num{163800}\left(\dfrac{t}{t_f}\right)^{11} +\num{108108}\left(\dfrac{t}{t_f}\right)^{10} - \num{40040}\left(\dfrac{t}{t_f}\right)^9 + \num{6435}\left(\dfrac{t}{t_f}\right)^8 \,
\end{align*}
with initial time $t_0=\SI{0}{\second}$ and final transition time~$t_{\rm f}=\SI{10}{\second}$. 
Note that $\sigma(t_0) = 0$, and $\sigma(t_{\rm f}) = 1$.
The initial rotation is~$y_0=\SI{0}{\radian\per\second}$ and the final rotation is $y_{\rm f} = 4  \pi \, \rm rad \, s^{-1} \approx \SI{12.56}{\radian\per\second} $, which corresponds to two revolutions per second. 

The experimental setup is controlled using the LabVIEW software package. The feedforward control is computed by solving the inverse model DAEs~(\ref{eq:servo-constraintsDAE}) in real-time using the implicit Euler scheme~(\ref{eq:implEul}). The strategy is implemented in LabVIEW and runs in real-time at~$\SI{1}{kHz}$, i.e., the control input in each time step is computed by solving one time step of the algorithm~(\ref{eq:implEul}) within the available time of $\Delta t = \SI{1}{\milli\second}$. The current implementation limits the control loop frequency to this value. In the following, the results presented for~$\SI{1}{kHz}$ are obtained for an online feedforward control implementation, while the results for~$\SI{2}{kHz}$ are computed with an offline feedforward control. The offline feedforward control is obtained by solving the same set of equations beforehand and accessing a lookup table of the solution during the control action. The feedback control law~(\ref{eq:uFC}) is also implemented in LabVIEW based on the measurements of the instantaneous values of the system output $y(t)$. The system output is obtained using angular encoders and filtering its data to obtain an angular velocity.

\subsection{Feedforward control}

We first validate the pure feedforward control strategy based on servo-constraints described in \Cref{sec:ffw}. The following results are obtained for a control loop frequency of~$\SI{1}{kHz}$. We compute the feedforward control input for the nominal model, which is defined by the parameters in \Cref{tab:param} and by setting the unknown friction force to $F_{\rm f,1}
=0$. In order to compensate modeling errors, the feedforward control input~$u_{\rm ffw}$ of the nominal model is adapted with two tuning factors $f_{\rm act}$ and $f_{\rm fric}$, such that the actuator input is
\begin{equation} \label{eq:ffwinput}
\begin{aligned}
u(t) &= f_{\rm act}\, u_{\rm ffw}(t) + f_{\rm fric}\,.
\end{aligned}
\end{equation}
Here, the tuning factor $f_{\rm fric}$ is supposed to compensate static Coulomb friction in the motion~$\dot \vp_1$ and the motor constant $f_{\rm act}$ is supposed to compensate the unidentified motor model which describes the relationship between the unit of the actuator input (electrical current) and the force input $u_{\rm ffw}$ determined by the inverse model.

\begin{figure}[h]
     \centering
     \begin{subfigure}[b]{0.42\textwidth}
         \centering
         \includegraphics[scale=0.85,page=2]{export_ffw_torin.pdf}
         \caption{System input.}
         \label{fig:ffw_u}
     \end{subfigure}
    \begin{subfigure}[b]{0.57\textwidth}
         \centering
         \includegraphics[scale=0.85]{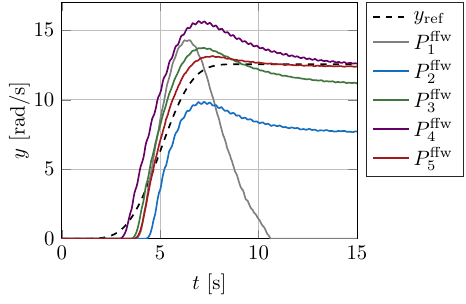}
         \caption{System output.}
         \label{fig:ffw_z}
     \end{subfigure}
        \caption{Real-time feedforward control at $\SI{1}{kHz}$.}
        \label{fig:ffw}
\end{figure}

The measurements are presented in \Cref{fig:ffw} for different tuning factors $f_{\rm act}$ and $f_{\rm fric}$. The various parameter sets for the tuning factors are listed in \Cref{tab:param_ffw}. The system input is shown in \Cref{fig:ffw_u}, while the measurement of the system output is shown in \Cref{fig:ffw_z}. Setting $f_{\rm fric}=0$ and therefore not compensating the friction in the system shows that the rotational velocity goes instantly back to zero due to friction in the actuator bearing. Therefore, different values for the Coulomb friction compensation are tested. Out of the taken measurements, the smallest tracking errors are achieved for the parameter set $P^{\rm ffw}_5$. These results show that the nominal model is sensitive to modeling errors. This is mainly due to the friction in the actuator. At the same time, the results validate the real-time application of the servo-constraints framework for a minimum phase system and demonstrate that accurate tracking is possible with a suitable feedforward controller. 

\begin{table}[h]
\begin{tabular}{ccc}
\hline 
& $f_{\rm act}$ & $f_{\rm fric}$ \\
\hline 
$P^{\rm ffw}_1$ & 0.3 & 0 \\[0.02cm]
$P^{\rm ffw}_2$ & 0.1 & 0.12 \\[0.02cm]
$P^{\rm ffw}_3$ & 0.1 & 0.15 \\[0.02cm]
$P^{\rm ffw}_4$ & 0.1 & 0.16 \\[0.02cm]
$P^{\rm ffw}_5$ & 0.08 & 0.16 \\
\hline
\end{tabular}
\caption{Tuning parameters for the feedforward controller.}
\label{tab:param_ffw}
\end{table}

\subsection{Feedback control}\label{sec:fb_control}
\begin{figure}[h]
     \centering
          \begin{subfigure}[b]{0.35\textwidth}
         \centering
         \includegraphics[scale=0.85,page=2]{export_fb_torin.pdf}
         \caption{System input.}
         \label{fig:fb_u}
     \end{subfigure}
     \hfill
     \begin{subfigure}[b]{0.63\textwidth}
         \centering
         \includegraphics[scale=0.85]{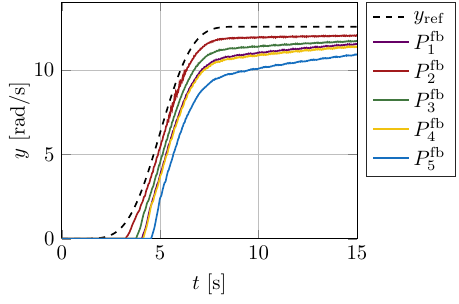}
         \caption{System output.}
         \label{fig:fb_z}
     \end{subfigure}
     \\[0.2cm]
     \begin{subfigure}[b]{0.35\textwidth}
         \centering
         \includegraphics[scale=0.85,page=3]{export_fb_torin.pdf}
         \caption{Tracking error for $P^{\rm fb}_2$ and $P^{\rm fb}_5$.}
         \label{fig:fb_e}
     \end{subfigure}
        \caption{Funnel feedback control at $\SI{2}{kHz}$.}
        \label{fig:fb}
\end{figure}

In order to accurately compensate for the unknown friction, we now apply the feedback control strategy described in \Cref{sec:fb}, without the feedforward control part. Since the funnel controller benefits from a high control loop frequency cf.~\cite{berger2020learning}, the following results are obtained for a frequency of $\SI{2}{kHz}$. We choose the exponential function
\begin{equation} \label{eq:psi}
\begin{aligned}
\psi(t) = s e^{-q t} + c
\end{aligned}
\end{equation}
to describe the performance funnel~$\psi$. The scalar parameters $s$, $q$ and $c$ change the size of the performance bound and act as design parameters. The chosen parameters are listed in \Cref{tab:param_fb} and the respective funnel functions are visualized in \Cref{fig:funneloverview}. The performance funnel $\psi_i$, which is plotted in the following figures, corresponds to the parameter set $P^{\rm fb}_i$. \Cref{fig:fb} shows the measurement results for the different parameter combinations~$P^{\rm fb}_i$. \Cref{fig:fb_u} shows the feedback control signal, while \Cref{fig:fb_z} shows the measured system output. All measurements have two properties in common: there is a time lag between the reference $y_{\rm ref}$ and the measured output~$y$ and there is a steady state error for $t>t_{\rm f}=\SI{10}{\second}$.
\Cref{fig:fb_e} shows the tracking error~$e=y - y_{\rm ref}$ and the performance funnel~$\psi$ of the best parameter set $P^{\rm fb}_2$ and the worst parameter set $P^{\rm fb}_5$. The results show that the tracking error stays close to the performance funnel $\psi$. For the tighter funnel, the tracking error comes very close to the boundary. This results in the chattering which is visible in the control input signal in \Cref{fig:fb_u}. High chattering motion can introduce undesirable loads on the mechanical parts of the actuator, and should therefore be avoided. 

\begin{table}[h]
\begin{tabular}{cccc}
\hline 
& $s$ & $q$ & $c$\\
\hline
$P^{\rm fb}_1$ &5 & 0.1 & 0.3 \\[0.04cm]
$P^{\rm fb}_2$ &1 & 0.1 & 0.5 \\[0.04cm]
$P^{\rm fb}_3$ &3 & 0.1 & 0.5 \\[0.04cm]
$P^{\rm fb}_4$ &5 & 0.1 & 0.5 \\[0.04cm]
$P^{\rm fb}_5$ &8 & 0.1 & 0.5 \\[0.04cm] \hdashline
$P^{\rm fb}_6$ &5 & 0.3 & 0.3 \\[0.04cm]
\hline
\end{tabular}
\caption{Parameter for the funnel function~\eqref{eq:psi}.  }
\label{tab:param_fb}
\end{table}

\begin{figure}[h]
     \centering
        \includegraphics[scale=0.85,page=1]{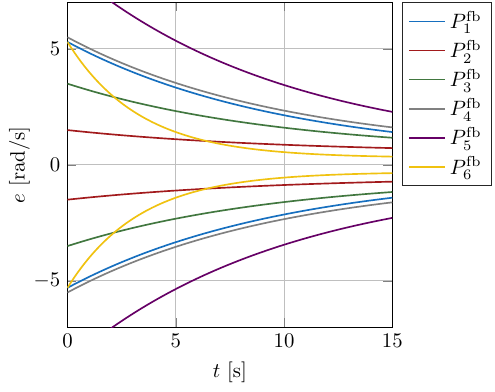}
        \caption{Visualization of the funnel parameters of \Cref{tab:param_fb}.}
        \label{fig:funneloverview}
\end{figure}

\begin{figure}[ht]
     \centering
         \begin{subfigure}[b]{0.35\textwidth}
         \centering
         \includegraphics[scale=0.85,page=2]{export_fb_and_ffw_torin5.pdf}
         \caption{System input.}
         \label{fig:comb1000_u}
     \end{subfigure}
     \hfill
     \begin{subfigure}[b]{0.64\textwidth}
         \centering
         \includegraphics[scale=0.85]{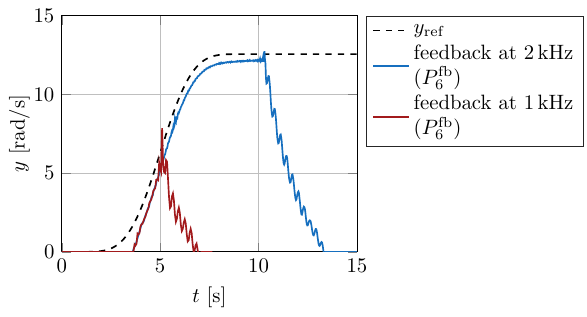}
         \caption{System output.}
         \label{fig:comb1000_z}
     \end{subfigure}
     \hfill
     \\[0.2cm]
    \begin{subfigure}[b]{0.35\textwidth}
         \centering
         \includegraphics[scale=0.85,page=3]{export_fb_and_ffw_torin5.pdf}
         \caption{Tracking error.}
         \label{fig:comb1000_e}
     \end{subfigure}
     \caption{Feedback control with parameter set $P^{\rm fb}_6$ at $\SI{1}{kHz}$ and $\SI{2}{kHz}$.}
     \label{fig:feedbackP6}
\end{figure}

In \Cref{tab:param_fb} another parameter set is listed, the application of which is not presented in \Cref{fig:fb}.
The reason of the absence of measurement results with~$P_{\rm fb}^6$ is that the feedback closed loop was unstable for both sample frequencies, because of the very aggressive funnel function. We plot the experimental results for the feedback controller with parameter set $P^{\rm fb}_6$ at $\SI{1}{kHz}$ and $\SI{2}{kHz}$ in \Cref{fig:feedbackP6}. The results show that both measurements yield the same tracking results as long as the controller is stable. However, the controller applied at $\SI{1}{kHz}$ becomes unstable after approximately $\SI{5}{\second}$, while the controller applied at $\SI{2}{kHz}$ becomes unstable after $\SI{10}{\second}$. Both controllers become unstable at some point, because the tracking error evolves too close to the performance bound and is pushed outside the limit by noise effects. However, as will be demonstrated in the following section, the combined controller is capable to achieve the control objective for parameter set~$P_{\rm fb}^6$.

\newpage
The combination of both, funnel feedback control and feedforward control, is presented in the following in order to minimize the drawbacks of the funnel feedback control. 

\subsection{Combination of both controller parts}
\begin{figure}[ht]
     \centering
         \begin{subfigure}[b]{0.35\textwidth}
         \centering
         \includegraphics[scale=0.85,page=2]{export_fb_and_ffw_torin3_rt.pdf}
         \caption{System input.}
         \label{fig:comb1000_u}
     \end{subfigure}
     \hfill
     \begin{subfigure}[b]{0.64\textwidth}
         \centering
         \includegraphics[scale=0.85]{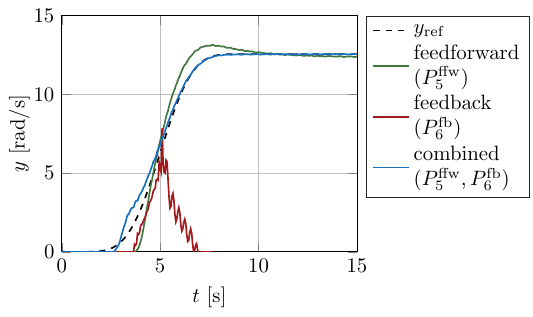}
         \caption{System output.}
         \label{fig:comb1000_z}
     \end{subfigure}
     \hfill
     \\[0.2cm]
    \begin{subfigure}[b]{0.35\textwidth}
         \centering
         \includegraphics[scale=0.85,page=3]{export_fb_and_ffw_torin3_rt.pdf}
         \caption{Tracking error.}
         \label{fig:comb1000_e}
     \end{subfigure}
     \caption{Combination of feedforward and feedback control at $\SI{1}{kHz}$.}
     \label{fig:comb1000}
\end{figure}

In the following, we present measurements for the combined control strategy and compare them to the results of the individual control strategies. 
Measurements are shown in \Cref{fig:comb1000} for the lower control loop frequency of $\SI{1}{kHz}$ with the funnel design parameters $P^{\rm fb}_6$ and the feedforward tuning parameters $P^{\rm ffw}_5$. We prefer to show these results over the results at the higher frequency of $\SI{2}{kHz}$, since the lower frequency is more realistic and applicable to a broader system class. Moreover, the feedforward solution can be computed in real-time for this case. The solution at the higher frequency of $\SI{2}{kHz}$ looks similar for the combined controller, but there are some differences in the application of the pure feedback controller, which are discussed in \Cref{sec:fb_control}. In the following figures, we compare the pure feedforward control, the pure feedback control and the combined strategy.

The results show that the pure feedback controller becomes unstable after a few seconds. \Cref{fig:comb1000_e} shows that the tracking error $e$ reaches the performance boundary~$\psi$ and then jumps out of the funnel between two time instances since the time discretization is not sufficiently small. However, applying the same parameter set $P^{\rm fb}_6$ for the funnel feedback controller and adding the feedforward control, results in a stable control signal. From this it can be concluded that the combined controller is amenable to operation at lower frequencies than pure feedback control. 

\subsection{Overall comparison}

After presenting individual measurements and pointing out different aspects of the analyzed control strategy, we now present an overall comparison of the control approaches. 
A qualitative comparison is shown in \Cref{tab:comparecontrollers}. The individual feedforward and feedback strategies have some disadvantages, respectively. The feedforward controller is, in the current implementation, real-time capable at the lower frequency of $\SI{1}{kHz}$. For the higher frequency, it can only be accessed via a pre-computed look-up table. Here, the real-time solution is preferred, since the feedforward controller can be adapted to varying trajectories or model parameters. On the other hand, the feedback controller benefits from higher control loop frequencies, since stricter performance bounds are possible. At the lower control loop frequency, the controller becomes unstable for realisitic performance bounds and the performance is only comparable to the feedforward controller. It therefore does not add any value at the lower frequency. The combined control strategy allows for even stricter performance bounds at both tested frequencies.

\begin{table}[h]
\begin{tabular}{|c|p{3cm}|p{3cm}|p{3cm}|}
\hline
 & \textbf{only feedforward} & \textbf{only feedback} & \textbf{combined controller} \\[0.1cm]
\hline
$\SI{2}{kHz}$ & \cellcolor{pink} access only via lookup table & \cellcolor{SpringGreen} stricter performance bounds possible compared to $\SI{1}{kHz}$ & \cellcolor{SpringGreen} robust application possible  \\[0.8cm]
\hline
$\SI{1}{kHz}$  & \cellcolor{SpringGreen} real-time computation possible & \cellcolor{pink}  unstable for realistic performance bounds &  \cellcolor{SpringGreen} robust application possible  \\[0.5cm]
\hline
\end{tabular}
\caption{Overview of qualitative controller properties.}
\label{tab:comparecontrollers}
\end{table}

\begin{figure}[h!]
     \centering
     \begin{subfigure}[b]{0.9\textwidth}
         \centering
         \includegraphics[scale=1]{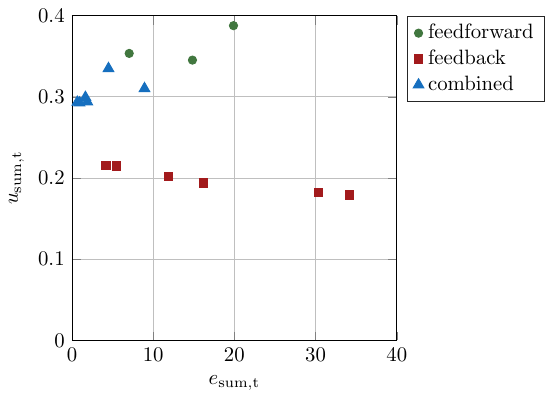}
         \caption{Transient regime.}
         \label{fig:performance_trans}
     \end{subfigure}   \\
     \begin{subfigure}[b]{0.9\textwidth}
         \centering
         \includegraphics[scale=1,page=2]{export_compareall.pdf}
         \caption{Stationary regime.}
         \label{fig:performance_stat}
     \end{subfigure}
        \caption{Performance measures in transient as well as stationary regime.}
        \label{fig:performance}
\end{figure}

Besides the qualitative comparison, we now introduce two performance measures in order to compare the overall performance of all performed experiments. First, we compare the performance in the transient regime between $t=\SI{0}{\second}$ and $t=\SI{10}{\second}$. We define the performance metrics
\begin{equation} \label{eq:errors_trans}
\begin{aligned}
u_{\rm sum}^{t_0,t_1} &:= \int_{t_0}^{t_1} u(\tau)^2 \, \text{d} \tau ,  \\
e_{\rm sum}^{t_0,t_1} &:= \int_{t_0}^{t_1} e(\tau)^2 \, \text{d} \tau 
\end{aligned}
\end{equation}
with the tracking error $e(t) = y(t) - y_{\rm ref}(t) $. 
To measure the total amount of input energy and the total tracking error in the transient regime we set $u_{\rm sum,t} := u_{\rm sum}^{0,t_f}$ and $e_{\rm sum, t} := e_{\rm sum,t }^{0,t_f}$, where $t_f = 10\, \rm s$ and the subscript~$\rm t$ indicates the transient regime. 
For integration of the measurement data we used the \textsc{Matlab} routine \textsc{trapz} with equidistant grid, since the sampling rate is fixed.
The transient performance comparison is shown in \Cref{fig:performance_trans}. The results show that the pure feedforward control signal uses comparatively large input energy while resulting in medium-sized tracking errors. The pure feedback controller uses less input energy with varying tracking errors. The combined control strategy uses more input energy compared to the feedback control strategy, but also results in smaller overall tracking errors, or to comparable tracking errors compared to the pure feedback controller. Thus, a direct comparison between the combined controller and the feedback controller is not possible based on these two chosen metrics. However, note that the pure feedback controller results in larger steady state errors which are not covered in this comparison of the transient regime. 
Therefore, we now compare the performance in the stationary regime between $t=\SI{10}{\second}$ and $t=\SI{15}{\second}$. 
To this end we consider the variance of the control input
\begin{equation} \label{eq:errors_stat_2}
\begin{aligned}
\sigma^2_{\rm sum,s}(u) &= \dfrac{1}{N} \sum_{i=1}^{N}(u(t_i) - \bar{u})^2 \quad \text{with} \quad
\bar  u = \dfrac{1}{N} \sum_{i=1}^{N} u(t_i) \,  \,
\end{aligned}
\end{equation}
where we identify $t_1 = t_f = 10 \, \rm s$ and $t_N = t_f + 5 \, \rm s = 15\, \rm s$.
The chattering of the control signal~$u$ is estimated by computing the variance~$\sigma^2_{\rm sum,s}$.
The total tracking error in the stationary regime is given by $e_{\rm sum,s} : = e_{\rm sum}^{10,15}$ according to~\eqref{eq:errors_trans} with the subscript s denoting the stationary regime.  

The stationary performance comparison is shown in \Cref{fig:performance_stat}. It is evident that the feedforward control strategy has no variance in the control signal while leading to large tracking errors $e_{\rm sum,s}$. In contrast, the feedback control signal has a comparatively large variance in the control input signal, while the stationary tracking errors $e_{\rm sum,s}$ vary in a large region. The combined control strategy exhibits the lowest variance in the control signal as well as the lowest tracking errors.
The clustering of two error values~$e_{\rm sum, t}$ each for the pure feedback in \Cref{fig:performance_trans} correspond to increasing offset values~$c$ of the funnel function~\eqref{eq:psi}, i.e., the smaller~$c$, the smaller $e_{\rm sum, t}$. Remarkably, the control action~$u_{\rm sum, t}$ is almost constant for different funnel functions, with slightly decreasing values of~$u_{\rm sum, t}$ for larger diameters of the funnel.
Note that some measurements with pure feedforward, which did not achieve the control objective may be not present in \Cref{fig:performance_trans} by choice of the image section.

This overall comparison supports the finding of the individual experimental results discussed above: combining funnel control with a feedforward controller can reduce chattering and enables smaller tracking errors by allowing smaller performance bounds. 

\section{Conclusion and outlook}\label{sec:conclusion}
In this contribution, we present experimental data for the combination of a feedforward controller based on servo-constraints and a feedback controller based on funnel control. A torsional oscillator serves as an experimental setup. 
Regarding the application of the method of servo-constraints, we validate its effectiveness for a real-time implementation on a minimum phase system. 
Regarding the combination of feedforward and funnel control, the experimental results support the theoretical work we obtained previously. The results show that this combination retains the advantages of both individual methods. In particular, we are able to apply funnel control at lower sample frequencies in the combined strategy compared to the application of pure funnel control. Moreover, we see a reduction of the peaky input signals, which are a known drawback of the funnel controller. Both of these advantages can lead to an easier application of funnel control in real-world systems.

Future work will focus on the implementation of the sample-and-hold feedback controller~\cite{lanza2023sampled} in the experimental setup.
Since this controller is designed to take only sampled system measurements, its application seems reasonable as well as promising.
Further, the feedforward controller can be improved during the control phase by implementing a learning scheme to improve the underlying model of the system.
Since the resulting controller is safeguarded by the funnel control component, this is safe learning in a real experimental setup.
A further extension of the combined controller described in the present article is the incorporation of a predictive controller, e.g., funnel MPC~\cite{denn2022funnel}, 
or its robust variant~\cite{berger2023robust}.

Moreover, as can be seen from the summarizing \Cref{tab:comparecontrollers}, applicability and robustness of pure feedforward/feedback as well as applicability of the combined controller strongly depends on the sampling frequency.
Therefore, in future research, we will conduct a series of experiments to systematically determine lower and upper limits for the sampling frequency to ensure safe operation of the combined controller.

\backmatter

\section*{Statements and Declarations}

\bmhead{Funding}
This work was funded by the Deutsche Forschungsgemeinschaft
(DFG, German Research Foundation) -- Project-IDs 362536361 (Thomas Berger, Svenja Drücker, Timo Reis and Robert Seifried), 396289190 (Robert Seifried, Svenja Drücker) and  471539468 (Lukas Lanza).

\bmhead{Author Contributions}
S. Drücker implemented and performed the experiments. L. Lanza provided the code for controller implementation. T. Berger, L. Lanza and T. Reis jointly provided the mathematical background for the control concept. The initial version of the manuscript was written by S. Drücker and L. Lanza with the following responsibilities: Section 1 and 2: jointly written, Section 3.1: S. Drücker, Section 3.2: L. Lanza, Section 4 and 5: S. Drücker, Section 6: jointly written. All
authors jointly conceived the presented ideas and analyses. All authors equally reviewed the manuscript. R. Seifried supervised the project and supported the experimental work.

\bmhead{Competing Interests} Non-financial interests: R. Seifried serves as an associate editor for the journal Multibody System Dynamics. The remaining authors have no conflicts of interest to declare that are relevant to the content of this article.









\bibliography{references}


\end{document}